\documentclass[12pt,a4paper,leqno]{article}

\usepackage{latexsym}
\usepackage{amssymb,exscale}
\usepackage[centertags]{amsmath}
\usepackage{amsthm}
\usepackage{graphicx}

\usepackage[dvips]{hyperref}

\numberwithin{equation}{section}

\theoremstyle{definition}
\newtheorem{example}{Example}
\newtheorem{theorem}{Theorem}
\swapnumbers
\newtheorem{remark}[equation]{Remark}

\renewcommand{\phi}{\varphi}

\newcommand{\I}{{\rm i}}
\newcommand{\D}{\mathrm{d}}
\newcommand{\E}{\mathrm{e}}

\renewcommand{\(}{\bigl(}
\renewcommand{\)}{\bigr)\vphantom{)}}

\newcommand{\ip}[2]{\langle#1,#2\rangle}

\newcommand{\ti}{\tilde}

\newcommand{\sgn}{\operatorname{sgn}}
\newcommand{\const}{\operatorname{const}}

\newcommand{\dist}{\operatorname{dist}}

\newcommand{\R}{\mathbb R}
\newcommand{\C}{\mathbb C}
\newcommand{\Z}{\mathbb Z}
\newcommand{\T}{\mathbb T}

\begin{document}

\title{Some extremal problems\\ related to Bell-type inequalities}

\author{Boris Tsirelson}

\date{}
\maketitle

\begin{abstract}
The best approximation by bounded product functions is calculated
for some very simple two-valued functions of two variables.
\end{abstract}

\section*{Introduction}
Here is a finite-dimensional extremal problem: given an $ n \times n
$\nobreakdash-matrix $ (a_{i,j})_{i,j} $ of numbers $ a_{i,j} = \pm1
$, maximize $ \sum_{i,j} a_{i,j} b_i c_j $ over all $ n
$\nobreakdash-vectors $ (b_i)_i $, $ (c_j)_j $ of numbers $ b_i
= \pm 1 $, $ c_j = \pm 1 $.

\begin{sloppypar}
The corresponding infinite-dimensional problem is: given a measurable
function of two variables $ f : X \times Y \to \{-1,1\} $, maximize $
\iint f(x,y) g(x) h(y) \, \D x \D y $ over all measurable functions $
g : X \to \{-1,1\} $, $ h : Y \to \{-1,1\} $. Here $ X,Y $ are given
measure spaces of finite measure. 
\end{sloppypar}

Equivalently, we seek the best $ L_2 $\nobreakdash-approximation of a
given function $ f(\cdot,\cdot) $ by factorizable functions $ g(\cdot)
h(\cdot) $ (all values being $ \pm 1 $).

More generally, we may consider measurable
vector\nobreakdash-functions $ g : X \to \R^d $, $ h : Y \to \R^d $
(for a given dimension $ d $) satisfying $ |g(x)| = 1 $, $ |h(y)| = 1
$ for all $ x,y $ ($ |\cdot| $ stands for the Euclidean norm). In this
case $ g(x) h(y) $ is interpreted as the inner product; that is, we
maximize
\[
\iint f(x,y) \ip{ g(x) }{ h(y) } \, \D x \D y \, .
\]
(Still, $ f : X \times Y \to \{-1,1\} $.) The case $ d=1 $ is just the
scalar case considered above. We may also use an
infinite\nobreakdash-dimensional separable Hilbert space $ H $ in
place of $ \R^d $ (the case $ d=\infty $).

About a relation to Bell-type inequalities see for example
\cite{Ts93}, especially (2.18), (2.20) and (2.24)--(2.25). The case $
d=1 $ is related to classical Bell-type inequalities; $ d=\infty $ ---
to quantum Bell-type inequalities in general; $ d=2 $ --- to quantum
Bell-type inequalities for a maximally entangled pair of qubits (the
singlet state of two spin\nobreakdash-$ 1/2 $ particles).

Recently some physicists \cite{phy} got especially interested in two
examples (see below) of such extremal problems. They conjectured the
optimal functions $ g,h $ (in both examples, for all $ d $) but did
not prove optimality of these functions. For $ d=1 $ they deduce their
claims from Inequality (7) of \cite{gi} (for $ n\to\infty $). However,
looking at \cite{gi} I did not find a proof of (7); it is just checked
by inspection for some small $ n $. My goal is to prove the optimality.

\begin{example}\label{ex1}
$ X = Y = (0,1) $;
\[
\begin{gathered}[m]\includegraphics{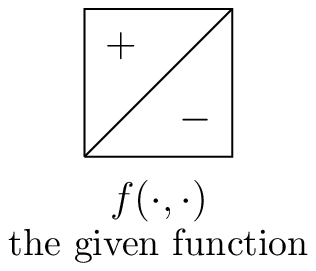}\end{gathered}
\quad
\begin{gathered}[m]
 f(x,y) = \begin{cases}
  1 &\text{if $ x<y $},\\
  -1 &\text{otherwise}.
 \end{cases}
\end{gathered}
\quad
\begin{gathered}[m]\includegraphics{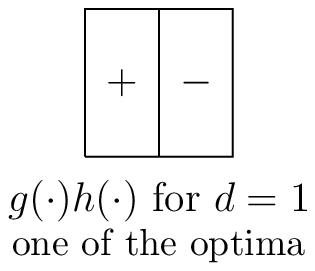}\end{gathered}
\]
\end{example}

\begin{example}\label{ex2}
$ X = Y = \T = \R/\Z $ (the circle of length $ 1 $);
\[
\begin{gathered}[m]\includegraphics{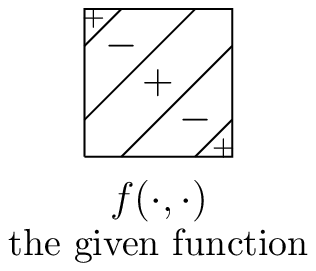}\end{gathered}
\begin{gathered}[m]
f(x,y) = \begin{cases}
 1 &\text{if $ \dist(x,y) < 0.25 $},\\
 -1 &\text{otherwise}.
\end{cases}
\end{gathered}
\begin{gathered}[m]\includegraphics{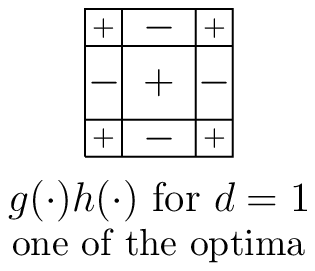}\end{gathered}
\]
\end{example}

The following result holds for Example \ref{ex1} as well as for
Example \ref{ex2}.

\begin{theorem}
\[
\max_{g,h} \iint f(x,y) \ip{ g(x) }{ h(y) } \, \D x \D y =
\begin{cases}
 \frac12 &\text{for $ d=1 $},\\
 \frac2\pi &\text{for $ 2 \le d \le \infty $}.
\end{cases}
\]
(The maximum is reached for every $ d $.)
\end{theorem}

\section[Some necessary conditions of optimality]{\raggedright Some
 necessary conditions of optimality}
\label{sect1}
Let functions $ g,h $ maximize $ I_f(g,h) = \iint f(x,y) \ip{ g(x) }{
h(y) } \, \D x \D y $ over all measurable $ g : X \to \R^d $, $ h : Y
\to \R^d $ (for a given dimension $ d $) satisfying $ |g(x)| = 1 $, $
|h(y)| = 1 $ for all $ x,y $. We introduce $ G : X \to \R^d $, $ H : Y
\to \R^d $ by
\begin{align*}
G(x) &= \int f(x,y) h(y) \, \D y \, , \\
H(y) &= \int f(x,y) g(x) \, \D x \, .
\end{align*}
Clearly, $ \int \ip{ g(x) }{ G(x) } \, \D x = I_f(g,h) = \int \ip{
H(y) }{ h(y) } \, \D y $. The optimality of $ g $ (for the given $ h
$) implies that $ g(x) = G(x) / |G(x)| $ for almost all $ x $ such
that $ G(x) \ne 0 $. Thus, $ I_f(g,h) = \int |G(x)| \, \D x
$. Similarly, $ h(y) = H(y) / |H(y)| $ for almost all $ y $ such
that $ H(y) \ne 0 $, and $ I_f(g,h) = \int |H(y)| \, \D y $.

We turn to Example \ref{ex1}: $ X = Y = (0,1) $ and $ f(x,y) =
\sgn(y-x) $ a.e. (almost everywhere). We have $ G(x) = -\int_0^x h(y)
\, \D y + \int_x^1 h(y) \, \D y $, therefore $ G $ is absolutely
continuous and $ G'(x) = -2h(x) $ a.e. Also, $ G(0) + G(1) = 0
$. Similarly, $ H'(y) = 2g(y) $ a.e., and $ H(0) + H(1) = 0 $. We get
a system of differential equations:
\begin{equation}\label{*}
\begin{aligned}
G'(x) &= -2 \frac{ H(x) }{ |H(x)| } \, , \\
H'(x) &= 2 \frac{ G(x) }{ |G(x)| } \, .
\end{aligned}
\end{equation}
Each equation holds almost everywhere, except for the points where its
right-hand side is undefined.

\begin{remark}
It can be deduced that $ |G(x)| + |H(x)| = \const $ and therefore
$ |G(x)| + |H(x)| = 2 I_f(g,h) $ for all $ x \in (0,1) $. However,
this fact will not be used.
\end{remark}

Now we turn to Example \ref{ex2}: $ X = Y = \T $ and $ f(x,y) = \sgn
(0.25-\dist(x,y)) $ a.e. We have $ G(x) = \int_{x-0.25}^{x+0.25} h(y)
\, \D y - \int_{x+0.25}^{x+0.75} h(y) \, \D y $, therefore $ G $ is
absolutely continuous and $ G'(x) = 2 h(x+0.25) - 2 h(x-0.25) $
a.e. (since $ x+0.75=x-0.25 $ in $ \T $). Also, $ G(x+0.5) = -G(x)
$. Similarly, $ H(x+0.5) = -H(x) $. We have $ G'(x) = 4 H(x+0.25) /
|H(x+0.25)| $ and $ H'(x) = 4 G(x+0.25) / |G(x+0.25)| $. Introducing $
\ti H(x) = H(x+0.25) $ we get a system of differential equations:
\begin{equation}\label{**}
\begin{aligned}
G'(x) &= 4 \frac{ \ti H(x) }{ |\ti H(x)| } \, , \\
\ti H'(x) &= -4 \frac{ G(x) }{ |G(x)| } \, .
\end{aligned}
\end{equation}
Each equation holds almost everywhere, except for the points where its
right-hand side is undefined.

\section[Dimension one]{\raggedright Dimension one}
\label{sect2}
Let $ d=1 $.

Differential equations \eqref{*} describe a dynamics on the plane $ \R^2
$; the integral curves evidently are the squares $ |G| + |H| = c $, $
c \in [0,\infty) $.
\[
\includegraphics{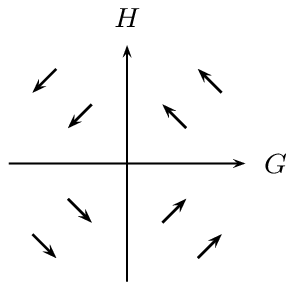}\qquad\includegraphics{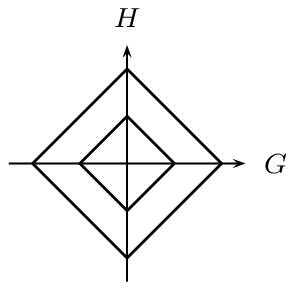}
\]
The solutions are periodic: $ G(x+2c) = G(x) $, $ H(x+2c) = H(x)
$. Also, $ G(x+c) = -G(x) $, $ H(x+c) = -H(x) $. The condition $
G(0)+G(1) = 0 $, $ H(0)+H(1) = 0 $ selects a sequence of solutions: $
c \in \{ 1, \frac13, \frac15, \dots \} \cup \{0\} $. We have to choose
$ c $ as to maximize $ I_f(g,h) $. Using the equality $ \int_0^1
|G(x)| \, \D x = I_f(g,h) = \int_0^1 |H(x)| \, \D x $ we get $
I_f(g,h) = 0.5 \int_0^1 \( |G(x)|+|H(x)| \) \, \D x = 0.5c $. The
maximizer is $ c=1 $ and the maximum is $ I_f(g,h) = 0.5 $.

The starting point $ \( G(0), H(0) \) $ can be chosen arbitrarily on
the square $ |G(0)| + |H(0)| = 1 $.
\[
 \includegraphics{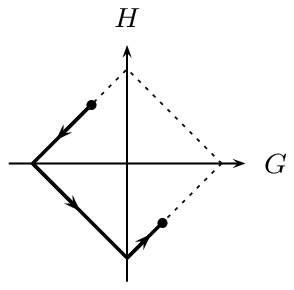}\qquad\includegraphics{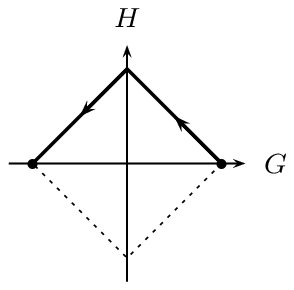}
\]
For instance, choosing $ G(0)=1 $, $ H(0)=0 $ we get $ H(\cdot) > 0 $
on $ (0,1) $, $ G(\cdot) > 0 $ on $ (0,0.5) $ and $ G(\cdot) < 0 $ on
$ (0.5,1) $. Thus, $ h(\cdot) = 1 $ on $ (0,1) $, $ g(\cdot) = 1 $ on
$ (0,0.5) $ and $ g(\cdot) = -1 $ on $ (0.5,1) $. This is the solution
shown on the picture (see Example \ref{ex1}).

Differential equations \eqref{**} are quite similar: $ |G|+|\ti H| = c
$; the condition $ G(x+0.5) = -G(x) $, $ \ti H(x+0.5) = -\ti H(x) $
selects $ c \in \{ 1, \frac13, \frac15, \dots \} \cup \{0\} $ again. The
maximizer is $ c=1 $ and the maximum is $ I_f(g,h) = 0.5 $.
\[
 \includegraphics{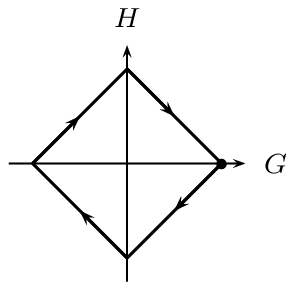}
\]
Choosing $ G(0)=1 $, $ \ti H(0)=0 $ we get $ G(\cdot) > 0 $ on $
(-0.25,0.25) $, $ G(\cdot) < 0 $ on $ (0.25,0.75) $, $ \ti H(\cdot) <
0 $ on $ (0,0.5) $ and $ \ti H(\cdot) > 0 $ on $ (-0.5,0) $. Thus, $
g(\cdot) = 1 $ on $ (-0.25,0.25) $, $ g(\cdot) = -1 $ on $ (0.25,0.75)
$, $ h(\cdot) = 1 $ on $ (-0.25,0.25) $ and $ h(\cdot) = -1 $ on $
(0.25,0.75) $. This is the solution shown on the picture (see Example
\ref{ex2}).

\section[Higher dimensions]{\raggedright Higher dimensions}
\label{sect3}
We start with Example \ref{ex2} and use Fourier series: $ G,H : \T \to
\R^d $,
\begin{gather*}
G(x) = \sum_k a_k \E^{2\pi\I kx} \, , \quad H(x) = \sum_k b_k
 \E^{2\pi\I kx} \, ; \\
a_k = \int_0^1 G(x) \E^{-2\pi\I kx} \, \D x  \, , \quad b_k = \int_0^1
 H(x) \E^{-2\pi\I kx} \, ;
\end{gather*}
$ a_k,b_k \in \C^d $ for $ k \in \Z $. We know that $ G(x+0.5) = -G(x)
$ and $ H(x+0.5) = -H(x) $, therefore $ a_k = 0 $, $ b_k = 0 $ for all
even $ k $. Especially, $ a_0 = 0 $, $ b_0 = 0 $. Taking into account
that
\begin{gather*}
\int_0^1 |G(x)|^2 \, \D x = \sum_k |a_k|^2 \, , \\
\int_0^1 |G'(x)|^2 \, \D x = (2\pi)^2 \sum_k k^2 |a_k|^2 \, ,
\end{gather*}
we get
\[
\int_0^1 |G'(x)|^2 \, \D x \ge 4\pi^2 \int_0^1 |G(x)|^2 \, \D x \, .
\]
According to \eqref{**}, $ |G'(x)| = 4 $ for almost all $ x $,
therefore $ \int_0^1 |G(x)|^2 \, \D x \le 4/\pi^2 $ and
\[
\int_0^1 |G(x)| \, \D x \le \frac2\pi \, .
\]
Using the equality $ I_f(g,h) = \int_0^1 |G(x)| \, \D x $ we get
\[
I_f(g,h) \le \frac2\pi \, .
\]
This bound can be reached already for $ d=2 $. Namely, we may take
two-dimensional vector-functions
\[
G(x) = \Big( \frac2\pi \cos 2\pi x, -\frac2\pi \sin 2\pi x \Big) \, ,
\quad H(x) = \Big( -\frac2\pi \sin 2\pi x, -\frac2\pi \cos 2\pi x
\Big)
\]
satisfying \eqref{**} and the conditions $ G(x+0.5) = -G(x) $, $
H(x+0.5) = -H(x) $. Then $ |G(x)| = 2/\pi $ for all $ x $, therefore $ 
I_f(g,h) = \frac2\pi $. More explicitly,
\[
g(x) = \( \cos 2\pi x, -\sin 2\pi x \) \, , \quad h(x) = \( -\sin 2\pi
x, -\cos 2\pi x \) \, .
\]

We turn to Example \ref{ex1}; here $ G,H : [0,1] \to \R^d $, $
G(0)+G(1) = 0 $, $ H(0)+H(1) = 0 $. We extend $ G,H $ to $ [0,2] $
letting
\[
G(1+x) = -G(x) \, , \quad H(1+x) = -H(x) \quad \text{for } x \in [0,1]
\]
and use Fourier series
\[
G(x) = \sum_k a_k \E^{\pi\I kx} \, , \quad H(x) = \sum_k b_k
 \E^{\pi\I kx}
\]
where $ a_k = 0 $, $ b_k = 0 $ for all even $ k $. We have
\begin{gather*}
\int_0^1 |G(x)|^2 \, \D x = \sum_k |a_k|^2 \, , \quad \int_0^1
 |G'(x)|^2 \, \D x = \pi^2 \sum_k k^2 |a_k|^2 \, , \\
4 = \int_0^1 |G'(x)|^2 \, \D x \ge \pi^2 \int_0^1 |G(x)|^2 \, \D x \,
, \quad \int_0^1 |G(x)| \, \D x \le \frac2\pi \, .
\end{gather*}
The bound is reached (in particular) for the
two\nobreakdash-dimensional vector\nobreakdash-functions
\[
G(x) = \Big( \frac2\pi \cos \pi x, \frac2\pi \sin \pi x \Big) \, ,
\quad H(x) = \Big( \frac2\pi \sin \pi x, -\frac2\pi \cos \pi x \Big)
\, ,
\]
which means
\[
g(x) = \( \cos \pi x, \sin \pi x \) \, , \quad h(x) = \( \sin \pi x,
-\cos \pi x \) \, .
\]

\bigskip
\filbreak
{
\small
\begin{sc}
\parindent=0pt\baselineskip=12pt
\parbox{4in}{
Boris Tsirelson\\
School of Mathematics\\
Tel Aviv University\\
Tel Aviv 69978, Israel
\smallskip
\par\quad\href{mailto:tsirel@post.tau.ac.il}{\tt
 mailto:tsirel@post.tau.ac.il}
\par\quad\href{http://www.tau.ac.il/~tsirel/}{\tt
 http://www.tau.ac.il/\textasciitilde tsirel/}
}

\end{sc}
}
\filbreak

\end{document}